\documentclass[a4paper, 10pt]{amsart}

\usepackage{etoolbox}
\makeatletter
\let\ams@starttoc\@starttoc
\makeatother
\usepackage[parfill]{parskip}
\makeatletter
\let\@starttoc\ams@starttoc
\patchcmd{\@starttoc}{\makeatletter}{\makeatletter\parskip\z@}{}{}
\makeatother

\usepackage{color,soul}
\definecolor{red}{rgb}{1,0,0}
\setulcolor{red}

\usepackage[colorlinks=true,linkcolor=blue,citecolor=blue,urlcolor=blue]{hyperref}
\usepackage{amsthm}
\usepackage[capitalise]{cleveref}
\usepackage{bookmark}

\usepackage{amssymb}
\usepackage{amsmath}
\usepackage{fancyhdr}
\usepackage{esint}
\usepackage{enumerate}

\usepackage{pictexwd,dcpic}

\swapnumbers
\newtheorem{lemma}{Lemma}[section]
\newtheorem{prop}[lemma]{Proposition}
\newtheorem{thm}[lemma]{Theorem}
\newtheorem{cor}[lemma]{Corollary}

\theoremstyle{definition}
\newtheorem{defn}[lemma]{Definition}
\newtheorem{example}[lemma]{Example}
\newtheorem{rem}[lemma]{Remark}
\newtheorem{ass}[lemma]{Assumption}

\numberwithin{equation}{section}

\renewcommand{\(}{\left(}
\renewcommand{\)}{\right)}

\renewcommand{\-}{\bar}

\newcommand{\cn}{\colon}

\newcommand{\R}{\mathbb{R}}

\renewcommand{\S}{\mathbb{S}}
\renewcommand{\H}{\mathbb{H}}

\renewcommand{\a}{\alpha}

\renewcommand{\d}{\delta}
\newcommand{\e}{\epsilon}
\renewcommand{\k}{\kappa}
\renewcommand{\l}{\lambda}

\newcommand{\vt}{\vartheta}

\newcommand{\s}{\sigma}

\newcommand{\G}{\Gamma}

\newcommand{\del}{\partial}

\newcommand{\fr}[2]{\frac{#1}{#2}}

\newcommand{\Thm}{\begin{thm}}
\newcommand{\eThm}{\end{thm}}
\newcommand{\Def}{\begin{defn}}
\newcommand{\eDef}{\end{defn}}
\newcommand{\Prop}{\begin{prop}}
\newcommand{\eProp}{\end{prop}}
\newcommand{\Rem}{\begin{rem}}
\newcommand{\eRem}{\end{rem}}
\newcommand{\Lem}{\begin{lemma}}
\newcommand{\eLem}{\end{lemma}}
\newcommand{\eq}{\begin{equation}}
\newcommand{\eeq}{\end{equation}}
\newcommand{\Ex}{\begin{example}}
\newcommand{\eEx}{\end{example}}
\newcommand{\pf}{\begin{proof}}
\newcommand{\epf}{\end{proof}}
\newcommand{\Cor}{\begin{cor}}
\newcommand{\eCor}{\end{cor}}
\newcommand{\Ass}{\begin{ass}}
\newcommand{\eAss}{\end{ass}}
\newcommand{\SAl}{\begin{align}\begin{split}}

\newcommand{\ra}{\rightarrow}
\newcommand{\hra}{\hookrightarrow}

\newcommand{\hp}{\hphantom}

\protected\def\ignorethis#1\endignorethis{}
\let\endignorethis\relax

\begin{document}

\title[Gradient estimates for ICF in hyperbolic space]{Gradient estimates for inverse curvature flows in hyperbolic space}
\author{Julian Scheuer}
\begin{abstract}
We prove gradient estimates for hypersurfaces in the hyperbolic space $\H^{n+1},$ expanding by negative powers of a certain class of homogeneous curvature functions. We obtain optimal gradient estimates for hypersurfaces evolving by certain powers $p>1$ of $F^{-1}$ and smooth convergence of the properly rescaled hypersurfaces. In particular, the full convergence result holds for the inverse Gauss curvature flow of surfaces without any further pinching condition besides convexity of the initial hypersurface.
\end{abstract}

\date{\today}
\thanks{This work is being supported by the DFG}
\address{Ruprecht-Karls-Universit{\"a}t, Institut f{\"u}r Angewandte Mathematik,
Im Neuenheimer Feld 294, 69120 Heidelberg, Germany}
\email{scheuer@math.uni-heidelberg.de}

\maketitle

\section{Introduction and the main result}

This short note is a direct improvement of \cite{Scheuer:06/2014}, in which we considered inverse curvature flows in the hyperbolic space $\H^{n+1},$ $n\geq 2,$ of the form
\SAl \label{floweq}		\dot{x}&=\fr{1}{F^{p}}\nu,\quad 0<p<\infty,\\
			x(0,M)&=M_{0},\end{split}\end{align}
where 
\eq x\cn [0,\infty)\times M\hra\H^{n+1}\eeq
is a family of embeddings, $F$ is a curvature function satisfying some additional properties specified later and $M_{0}$ is a suitable initial hypersurface. We showed that under certain conditions on $F$ and $M_{0}$ this flow exists for all time and converges to a well-defined smooth function after rescaling. In case $p>1$ we had to impose quite restrictive conditions on $F$ and $M_{0},$ namely $F$ was supposed to vanish on the boundary of 
\eq\G_{+}=\{(\k_{i})\in\R^{n}\cn \k_{i}>0\quad\forall i=1,\dots,n\}\eeq
and $M_{0}$ had to satisfy a pinching condition on the oscillation, cf. \cite[Thm. 1.2, (2)]{Scheuer:06/2014}. This condition was designed to ensure that the gradient of the initial convex hypersurface $M_{0}$ was small. Then it was possible to show that it remained small. The unsatisfactory thing was that we were not able to treat powers of the Gaussian curvature that were larger than $n^{-1}$ without this restriction. The aim of this note is to provide gradient estimates which allow to treat \eqref{floweq} for some powers $1<p\leq p_{0},$ where $p_{0}$ will depend on the curvature function. In case of the Gaussian curvature, $F=nK^{\fr 1n},$ we will obtain $p_{0}=\fr{n}{n-1}.$ This includes the full convergence result in case of the inverse Gauss curvature flow for surfaces with arbitrary convex initial hypersurface, a problem which to our knowledge has not been considered in the literature before. The class of curvature functions we are able to prove the main result for is closely related to the class $(K^{*})$ defined in \cite[Def. 2.2.15]{Gerhardt:/2006}. We will come back to this issue in some more detail later.

Before we state the main result, let us formulate the assumptions we have to impose on the curvature function. We have to add one assumption compared to \cite[Thm. 1.2]{Scheuer:06/2014}.

\Ass\label{Ass}
Let 
\eq\G_{+}=\{(\k_{i})\in\R^{n}\cn \k_{i}>0\quad\forall i=1,\dots,n\}.\eeq
Assume $F\in C^{\infty}(\G_{+})\cap C^{0}(\-{\G}_{+})$ to be a
\begin{itemize}
\item{symmetric}
\item{monotone}
\item{homogeneous of degree $1$} 
\item{concave}
\end{itemize}
 curvature function satisfying
\eq F_{|\del\G_{+}}=0,\eeq
 the normalization
\eq F(1,\dots,1)=n\eeq
and
\eq\label{Ass4} \fr{\del F}{\del \k_{i}}\k_{i}\geq \e_{0}F\quad\forall i=1,\dots,n\eeq
for some $0<\e_{0}(F)\leq \fr 1n.$
\eAss

We will prove the following result.

\Thm\label{Main}
Let $n\geq 2$ and 
\eq x_{0}\cn M\hra \H^{n+1}\eeq
be the smooth embedding of a closed and strictly convex hypersurface $M_{0}.$ Let $F$ satisfy \cref{Ass} and
\eq 1<p\leq p_{0}:=\fr{1}{1-\e_{0}}.\eeq
Then the following statements hold.
\begin{enumerate}[(i)]
\item{There exists a unique global solution
\eq x\cn [0,\infty)\times M\hra\H^{n+1}\eeq
of the curvature flow equation
\SAl	\dot{x}&=\fr{1}{F^{p}}\nu\\
		x(0,M)&=M_{0},\end{split}\end{align}
		where $\nu$ is the outward normal to the flow hypersurfaces $M_{t}=x(t,M)$ and $F$ is evaluated at the principal curvatures of $M_{t}.$}
\item{Representing $M_{0}$ as a graph in geodesic polar coordinates,
\eq M_{0}=\{(u(0,x^{i}),x^{i}(0,\xi))\cn \xi\in M\},\eeq
where $u=x^{0}$ describes the radial distance to a point inside the convex body enclosed by $M_{0},$ we obtain that all flow hypersurfaces $M_{t}$ have a similar representation by a scalar function 
\eq u\cn [0,\infty)\times \S^{n}\ra \R.\eeq
The rescaled hypersurfaces $\~M_{t}$ which are described via the scalar function
\eq\label{RescFlow} \~u=u-\fr{t}{n^{p}}\eeq
converge to a well-defined hypersurface in $C^{\infty}.$}
\item{In case 
\eq F=nK^{\fr 1n},\eeq
where $K$ is the Gaussian curvature, we have
\eq p_{0}=\fr{n}{n-1}.\eeq}
\end{enumerate}
\eThm

\Rem
Note that the convergence of the functions \eqref{RescFlow} can not be improved in general. A counterexample was recently given by Hung and Wang in case of the inverse mean curvature flow, $F=H$ and $p=1,$ cf. \cite{HungWang:06/2014}. Also in case $p>1$ we do not expect that the convergence can be improved, since the evolution equation of the traceless second fundamental form
\eq \|\mathring{A}\|^{2}=\|A\|^{2}-\fr 1n H^{2}\eeq
has the same structure in the first order terms as in case $p=1.$
\eRem

\Rem
Note that the $n$-th root of the Gaussian curvature
\eq F=nK^{\fr 1n}\eeq
fulfills \cref{Ass} with $\e_{0}=\fr 1n.$ To see that the set of curvature functions satisfying \cref{Ass} contains considerably more functions the reader is referred to \cite[Prop. 2.2.18]{Gerhardt:/2006}. In case $n=2$ this class contains, up to normalization and smoothness, all curvature function of class $(K^{*})$ which are homogeneous of degree $1.$ For a definition of the class $(K^{*})$ compare \cite[Sec. 2.2]{Gerhardt:/2006}. For arbitrary $n\geq 2$ it contains those $1$-homogeneous functions which can be written as
\eq F=GK^{a}, a>0,\eeq
where $G$ shares all the properties of a function belonging to class $(K)$ besides vanishing on $\del \G_{+}.$  

Furthermore note that the value $\e_{0}$ in \eqref{Ass4} can not be larger than $\fr 1n,$ due to the homogeneity, which implies
\eq F^{i}\k_{i}=F.\eeq
\eRem

\section{Notation}
Let us refresh some notation already used in \cite{Scheuer:06/2014}. We consider hypersurfaces in $\H^{n+1},$ the metric of which reads in geodesic polar coordinates 
\eq d\-s^{2}=dr^{2}+\vt^{2}(r)\s_{ij}dx^{i}dx^{j},\eeq
where 
\eq \vt(r)=\sinh(r),\eeq
$r=x^{0}$ denotes the geodesic distance to some point $q\in\H^{n+1}$ and $\s_{ij}$ is the round metric of the $n$-sphere. In such coordinates, let the hypersurface $M$ be given as a graph over a geodesic sphere $\S^{n}$ around $q,$
\eq M=\{(u(x),x^{i})\cn (x^{i})\in\S^{n}\}.\eeq
Then the induced metric of $M$ is given by
\eq g_{ij}=u_{i}u_{j}+\vt^{2}(u)\s_{ij}.\eeq
Let
\eq v^{2}=1+\vt^{-2}\s^{ij}u_{i}u_{j}\equiv 1+|Du|^{2},\eeq
 then the outward unit normal is given by
 \eq (\nu^{\a})=v^{-1}(1,-\vt^{-2}\s^{ij}u_{j}),\eeq
 $\a=0,\dots,n$ and $i=1,\dots,n.$
 For curvature functions depending on the second fundamental form and the metric,
 \eq F=F(h_{kl},g_{kl})\eeq
 we let
 \eq F^{kl}=\fr{\del F}{\del h_{kl}}.\eeq
 
 For functions $f$ defined on a manifold $M,$ lower indices indicate covariant differentiation with respect to the induced metric of $M,$ e.g. $u_{i}$ or $v_{ij}.$
 
  A dot over a function or a tensor always indicates time derivation, e.g.
 \eq \dot{v}=\fr{d}{dt}v,\eeq
 whereas a prime always denotes derivation with respect to a direct argument. For example, if $f=f(u),$ then
 \eq f'=\fr{d}{du}f\eeq
 and
 \eq \dot{f}=f'\dot{u}.\eeq
 Note that this notation partially deviates from those in \cite{Gerhardt:01/1996} and \cite{Gerhardt:/2006}.
 
 \section{Rough outline of the proof}\label{Outline}
 Let us shortly explain, which ingredient of our proof is the crucial one compared to the proof in \cite{Scheuer:06/2014}. In \cite[Sec. 3]{Scheuer:06/2014} we proved the longtime existence of the solution to \eqref{floweq} by successively proving the existence of spherical barriers and of bounds on $v$ and the second fundamental form. This part of the proof is not affected by our relaxed assumptions at all, compare \cite[Thm. 1.2 (1)]{Scheuer:06/2014}. Thus we have a longtime graph representation of the flow hypersurfaces,
 \eq M_{t}=\{(u(t,x^{i}),x^{i}(t,\xi))\cn \xi\in M\}.\eeq
 
 In order to derive convergence of the flow, we started the decay estimates by proving that the quantity $v$ converges to $1$ exponentially fast. In case $p>1$ this was only possible if $v$ was sufficiently small initially, due to the positively signed second order term in \cite[equ. (4.3)]{Scheuer:06/2014}. As we will show in \cref{Proof}, this restriction is not necessary in the situations described in \cref{Main}. It turned out to be possible to further exploit the first negatively signed term in \cite[equ. (4.3)]{Scheuer:06/2014} by using a better suited arrangement of the terms involved. We will be able to exploit this extra term with the help of the additional assumptions in \cref{Ass} to derive exponential decay of $v-1.$
 
In order to obtain the optimal gradient estimate 
\eq v-1\leq ce^{-\fr{2t}{n^{p}}}\eeq
and to conclude the final convergence result we observe that the further arguments in \cite{Scheuer:06/2014} do not depend on the pinching restriction at all and the proofs of the further decay estimates apply literally, once one has proven decay of $v-1.$

\section{Proof of the main theorem}\label{Proof}

As mentioned in \cref{Outline} it suffices to prove decay of $v-1.$ The following proposition holds.

\Prop
Let $x$ be the solution of \eqref{floweq} under \cref{Ass} with strictly convex initial hypersurface $M_{0}$ and
\eq 1<p\leq \fr{1}{1-\e_{0}}.\eeq
Let $u$ be a corresponding graph representation of the $M_{t}.$
Then the quantity $v^{2}=1+|Du|^{2}$ satisfies the decay estimate
\eq v-1\leq ce^{-\l t},\eeq
with suitable positive constants $c, \l$ which depend on $n,$ $M_{0},$ $p$ and $F.$
\eProp

\pf
The flow \eqref{floweq} is of the form
\eq \dot{x}=-\Phi(F)\nu,\eeq
where 
\eq \Phi(r)=-r^{-p}.\eeq

According to \cite[(5.28)]{Gerhardt:01/1996} for such a flow we have
\SAl	\label{Grad1}	\dot{v}-\Phi'F^{ij}v_{ij}&=-\Phi'F^{ij}h_{ik}h^{k}_{j}v-2v^{-1}\Phi'F^{ij}v_{i}v_{j}+2\Phi'F^{ij}v_{i}u_{j}\fr{\-H}{n}\\
						&\hp{=}-\Phi'F^{ij}g_{ij}\fr{\-H^{2}}{n^{2}}v-\Phi'F^{ij}u_{i}u_{j}\fr{{\-H'}}{n}v+\fr{\-H}{n}(\Phi-\Phi'F)|Du|^{2}\\
						&\hp{=}+2\Phi'F\fr{\-H}{n}v^{2}	\\
						&=-\Phi'F^{ij}\(h_{ik}h^{k}_{j}-2\fr{\-H}{n}h_{ij}+\fr{\-H^{2}}{n^{2}}g_{ij}\)v+\(1-\fr{1}{p}\)\Phi'F\fr{\-H}{n}v^{2}\\
						&\hp{=}-2\Phi'F\fr{\-H}{n}v+\(1+\fr{1}{p}\)\Phi'F\fr{\-H}{n}-\Phi'F^{ij}u_{i}u_{j}\fr{{\-H'}}{n}v\\
						&\hp{=}-2v^{-1}\Phi'F^{ij}v_{i}v_{j}+2\Phi'F^{ij}v_{i}u_{j}\fr{\-H}{n}.	\end{split}\end{align}
						
Define 
\eq w=(v-1)e^{\l t},\quad \l>0,\eeq
and suppose for $0<T<\infty$ that
\eq \sup\limits_{[0,T]\times M}w=w(t_{0},\xi_{0})\geq 1,\quad t_{0}>0.\eeq

Due to \cite[(5.29)]{Gerhardt:01/1996} we have at $(t_{0},\xi_{0})$
\eq 0=v_{i}=-v^{2}h^{k}_{i}u_{k}+v\fr{\-H}{n}u_{i}.\eeq
Choose Riemannian normal coordinates around $(t_{0},\xi_{0})$, in which
\eq g_{ij}=\d_{ij},\quad h_{ij}=\k_{i}\d_{ij}.\eeq
Since $v(t_{0},\xi_{0})>1,$ there exists an index $j\in\{1,\dots,n\},$ such that $u_{j}\neq 0$ and thus
\eq \k_{j}=v^{-1}\fr{\-H}{n}.\eeq

Noting that in the hyperbolic space the mean curvature $\-H$ of the slices $\{x^{0}=u\}$ satisfies
\eq \fr{\-H}{n}=\coth u,\eeq
that
\eq \fr{\del F}{\del \k_{i}}\geq \e_{0}\fr{F}{\k_{i}}\quad\forall 1\leq i\leq n\eeq
and that in our present coordinate system we have
\eq F^{ii}=\fr{\del F}{\del \k_{i}},\eeq
cf. \cite[Lemma 2.1.9]{Gerhardt:/2006},
we obtain from \eqref{Grad1} at the point $(t_{0},\xi_{0})$ that 
\SAl		0&\leq\l w-\Phi'F^{jj}\fr{\-H^{2}}{n^{2}}\fr{(v-1)^{2}}{v}e^{\l t_{0}}+\(1-\fr{1}{p}\)\Phi'F\fr{\-H}{n}v^{2}e^{\l t_{0}}\\
						&\hp{=}-2\Phi'F\fr{\-H}{n}ve^{\l t_{0}}+\(1+\fr{1}{p}\)\Phi'F\fr{\-H}{n}e^{\l t_{0}}-\Phi'F^{ij}u_{i}u_{j}\(1-\fr{\-H^{2}}{n^{2}}\)ve^{\l t_{0}}\\
						&\leq \l w-\e_{0}\Phi'F\fr{\-H}{n}(v-1)^{2}e^{\l t_{0}}+\(1-\fr{1}{p}\)\Phi'F\fr{\-H}{n}v^{2}e^{\l t_{0}}\\
						&\hp{=}-2\Phi'F\fr{\-H}{n}ve^{\l t_{0}}+\(1+\fr{1}{p}\)\Phi'F\fr{\-H}{n}e^{\l t_{0}}-\Phi'F^{ij}u_{i}u_{j}\(1-\fr{\-H^{2}}{n^{2}}\)ve^{\l t_{0}}\\
						&=\l w+\(1-\fr 1p-\e_{0}\)\Phi'F \fr{\-H}{n}v^{2}e^{\l t_{0}}-\(2-2\e_{0}\)\Phi'F\fr{\-H}{n}ve^{\l t_{0}}\\
						&\hp{=}+\(1+\fr{1}{p}-\e_{0}\)\Phi'F\fr{\-H}{n}e^{\l t_{0}}+ce^{\(\l-\fr{2}{n^{p}}\)t_{0}}\\
						&\leq \l w-\(1-\e_{0}\)\Phi'F\fr{\-H}{n}w+\(1-\fr{1}{p}-\e_{0}\)\Phi'F\fr{\-H}{n}vw+ce^{\(\l-\fr{2}{n^{p}}\)t_{0}}\\
						&<0\end{split}\end{align}
for small $\l$ and large $t_{0},$ where we also used 
\eq F+F^{-1}+v\leq c,\eeq
\eq \fr{\-H}{n}-1\leq ce^{-\fr{2t}{n^{p}}}\eeq
and that the $\k_{i}$ range in a compact set of $\G_{+}.$ We derived those facts in \cite[Cor. 3.6, Lemma 3.7, Prop. 3.10 and Prop. 3.11]{Scheuer:06/2014}. This is a contradiction and thus $w$ is bounded.
\epf

\Rem
In case of the Gaussian curvature
\eq F=nK^{\fr 1n}\eeq
we obtain 
\eq \fr{1}{1-\e_{0}}=\fr{n}{n-1},\eeq
which proves part (iii) of \cref{Main}.
\eRem

\bibliographystyle{/Users/J_Mac/Documents/Uni/TexTemplates/hamsplain}
\bibliography{/Users/J_Mac/Documents/Uni/TexTemplates/Bibliography}
\end{document}